\documentclass[a4paper,11pt]{article}
\usepackage[latin1]{inputenc}
\usepackage[T1]{fontenc}
\usepackage[english]{babel}
\usepackage{eurosym}
\usepackage[pdftex,usenames,dvipsnames]{color}
\usepackage[pdftex]{graphicx}

\usepackage{hyperref}

\usepackage{caption2}
\usepackage{amsmath}
\usepackage{amssymb}
\usepackage{fancybox}
\usepackage{amscd}
\usepackage{amsthm}

\theoremstyle{plain}
\newtheorem{Thm}{Theorem}

\theoremstyle{definition}

\newtheorem{Remark}{Remark}
\theoremstyle{remark}

\numberwithin{equation}{section}

\title{
Indices of
   the iterates of \({\Bbb R}^3\)-homeomorphisms at Lyapunov stable fixed points}

\author{Francisco R. Ruiz del Portal and José M. Salazar \thanks{ The authors have been supported by MEC,
MTM 2006-0825. \newline 2000 {\em Mathematics Subject
Classification}: 37C25, 37B30, 54H25.
}}

\begin{document}

\maketitle

\begin{abstract}
Given any positive sequence \(\{c_n\}_{n \in {\Bbb N}}\), we
construct orientation preserving homeomorphisms \(f:{\Bbb R}^3
\rightarrow  {\Bbb R}^3\) such that \(Fix(f)=Per(f)=\{ 0\}\),
\(0\) is Lyapunov stable and \( \limsup \frac{|i(f^m, 0)|}{c_m}=
\infty\). We will use our results to discuss and to point out some
strong differences with respect to the computation and behavior of
the sequences of the indices of planar homeomorphisms.
\end{abstract}


\maketitle

\centerline{{\bf 1. Introduction.}}

\medskip

The  computation of the sequence of the indices, or the sequence
of Lefschetz numbers, of the iterates of a map is an important and
non-trivial problem.

When a fixed point is an isolated invariant set of an orientation
preserving planar homeomorphism, the problem of the computation of
the indices of its iterates was solved by Le Calvez and Yoccoz,
(\cite{CY} and \cite{CY1}) and, by the authors, in the orientation
reversing case (\cite{RS}). Later Le Calvez solved the general
problem in the orientation preserving case using the
Carath\'eodory's theory of prime ends (\cite{LC}) and the authors,
in \cite{RS1}, the general case for orientation reversing planar
homeomorphisms.

For orientation preserving planar homeomorphisms there are
integers \(r\) and \(q\) such that the sequence of indices  is as
follows:

\[ i(f^{k}, p)=\left\{\begin{array}{ll}
1-rq & \text{ if } k \in r{\Bbb Z} \\
1 & \text{ if } k \notin r{\Bbb Z} \end{array} \right.
\]

If the problem in the plane resulted to be hard, the analogous
problem in \({\Bbb R}^3\) seems to be strongly non-trivial because
of the different dynamical pathologies that can appear. For
instance, while Lyapunov stable isolated fixed points of planar
homeomorphisms have always index =1 i.e. the Euler characteristic
of a disc, (\cite{DO} and \cite{RP}), for \({\Bbb R}^n\)-vector
fields
 with \(n \geq 3\), Bonatti and Villadelprat in \cite{BV} proved that the index
of stable, even in the past and in the future, isolated rest
points can be any integer (see also the paper of Erle, \cite{E}).

There are not many known results about the behavior of the
sequences  of fixed point indices of homeomorphisms in dimension
3. For instance it is well known that the  sequence must follow
Dold's necessary conditions (\cite{Do}). Shub and Sullivan proved
that for \(C^1\)-maps (no necessarily injective) the sequence is
bounded. Later, Chow, Mallet-Paret and Yorke (\cite{CMPY}) gave
bounds about the form of the sequence of indices in terms of the
spectrum of the derivative \(Df(p)\). Babenko and Bogatyi
(\cite{BB}) proved that these bounds are sharp in dimension 2 and
in a recent paper Graff and Nowak-Przygodzki have proved
(\cite{GN}) that for \(C^1\)-maps the sequence of fixed point
indices follows one among exactly seven different periodic
patterns.

More recently, the authors, see \cite{RS2}, have solved completely
the problem for \({\Bbb R}^3\)-homeomorphisms belonging to a
special class. A class that is quite natural to study because the
corresponding family in \({\Bbb R}^2\) is  the set of all planar
homeomorphisms such that \(\{p\}\) is an isolated invariant set.

Let \(U \subset {\Bbb R}^3\) be an open subset and let \(\cal{B}\)
be  the set of all homeomorphisms \(f: U \subset {\Bbb R}^3
\rightarrow f(U) \subset {\Bbb R}^3\) such that there exists a
closed 3-dimensional ball, \(N\), with the following properties:

a)  \(N\) is an isolating block such that the maximal invariant
set contained in \(N\),  \(Inv(N,f)\), is \( \{p\}\),

b) \(\partial N\) is a locally flat 2-sphere and

c) the component of \(f(N) \cap N\) containing \(p\) is also a
closed ball.

In \cite{RS2} it is shown that for every \(f\in \cal{B}\) the
sequence \(\{i(f^k,p)\}_{k \in {\Bbb N}}\), is periodic.
Conversely, for any periodic sequence of integers \(\{r_k\}_{k \in
{\Bbb N}}\) satisfying Dold's necessary congruences, there exists
an orientation preserving homeomorphism \(f \in \cal{B}\) such
that \(i(f^k,p)=r_k\) for every \(k \in {\Bbb N}\).

Shub and Sullivan in \cite{SS} also conjectured that for every
\(C^1\)-map, \(h\), defined in a compact manifold,

\[\limsup \frac{log(|Per^m(h)|)}{m} \geq \limsup
\frac{log(|\Lambda(h^m)|)}{m},\]

\noindent where \(\Lambda\) denotes the Lefschetz number.
Obviously every homeomorphism of the \(n\)-sphere, \(S^n\),
satisfies the above inequality because the sequence of the
Lefschetz numbers of its iterates is constant if it preserves
orientation.

A more general and slightly different version of the above problem
is whether

\[ \limsup \frac{log(|Per^m(h)|)}{m} \geq \limsup \frac {log(\sum_{p \in
Per^m(h)}|i(h^m,p)|)}{m}.\]

It is well known that there are examples of non-injective
continuous maps for which both previous inequalities fail to be
true (\cite{SS}).

In this paper we see that for \(S^3\)-homeomorphisms the answer to
this second problem is negative. We will show that, if for any
neighborhood \(N\) of the origin, \(Inv(N,f) \cap
\partial(N) \ne \emptyset\) then, the sequence of indices of the iterates
of \(f\) may be unbounded. In other words, if in the conjecture of
Shub and Sullivan we replace a manifold by a bounded open subset
of \({\Bbb R}^3\), a \(C^1\)-map by a homeomorphism and the
Lefschetz numbers by the fixed point indices, the answer is
negative even for stable fixed points. More precisely we shall
prove the following theorems that solve, in the negative, Problem
2.3.1 of \cite{Va}.

\begin{Thm}
For each positive sequence \(\{c_n\}_{n \in {\Bbb N}}\), there
exist orientation preserving \( {\Bbb R}^3\)-homeomorphisms,
\(f\), such that \(Fix(f)=Per(f)=\{0\}\) and \( \limsup
\frac{|i(f^m, 0)|}{c_m}= \infty\). Moreover, if \(B \subset {\Bbb
R}^3\) is any closed ball centered in the origin, \(f(B) \cap B\)
is a topological ball, \(Inv(B,f)\) is the closed 2-disc \(B \cap
\{z=0\}\) and \(f\) is limit of a sequence of homeomorphisms
\(\{f_m\}_{m \in {\Bbb N}}\) such that \(Inv(B,f_m)=\{ 0\}\) for
every \(m \in {\Bbb N}\) and, for every \(n \in {\mathbb N}\),
there exist \(m_0\) such that \( i(f^n, 0)= i((f_m)^n,0)\) for
every \(m \geq m_0\).
\end{Thm}

\begin{Thm}
For each positive sequence \(\{c_n\}_{n \in {\Bbb N}}\), there
exist orientation preserving \( {\Bbb R}^3\)-homeomorphisms,
\(h\), such that \(Fix(h)=Per(h)=\{0\}\), \(0\) is Lyapunov stable
and \( \limsup \frac{|i(h^m, 0)|}{c_m}= \infty\). In particular,
there are \({\Bbb R}^3\)-homeomorphisms, such that
\(Fix(h)=Per(h)=\{ 0\}\), \(0\) is Lyapunov stable and \( \limsup
\frac{log(|i(h^m, 0)|)}{m}= \infty\).
\end{Thm}

The techniques used for the computation of the indices are valid
for both, orientation preserving and orientation reversing
homeomorphisms.

\medskip

If \(X\) is a compact ANR (absolute neighborhood retract for
metric spaces), \(i_X(f, p)\) will denote, if it is well defined,
the fixed point index of \(f\) in a small enough neighborhood of
\(p\). When the indices are computed in the Euclidean space we
shall write just \(i(f,p)\).

The reader is referred to the text of \cite{Br2}, \cite{Do},
\cite{Nu} and the recent book of Jezierski and Marzantowicz,
\cite{JM},  for information about the fixed point index theory.
The last one is also appropriated to find in a unified way the
results of \cite{BB}, \cite{CMPY} and  \cite{SS} we mentioned
above.

\medskip
\medskip

\centerline{{\bf 2. Preliminary definitions and some basic
examples.}}

\medskip

Given \(A \subset B \subset N\), \(cl(A)\), \(cl_B(A)\),
\(int(A)\), \(int_B(A)\), \(\partial A\) and \(\partial_B A\) will
denote the closure of \(A\), the closure of \(A\) in \(B\), the
interior of \(A\), the interior of \(A\) in  \(B\), the boundary
of \(A\) and the boundary of \(A\) in \(B\) respectively.

Let \(U \subset X\) be an open set. By a {\em (local)
semidynamical system} we mean a local homeomorphism \(f:U
\rightarrow X\). The {\em invariant part} of \(N\), \(Inv(N, f)\),
is defined as the set of all \(x \in N\) such that there is a full
orbit \(\gamma\) with \(x \in \gamma \subset N\).

A compact set \(S\subset X\) is {\em invariant} if \(f(S)=S\). A
compact invariant set \(S\) is {\em isolated with respect to}
\(f\) if there exists a compact neighborhood \(N\) of \(S\) such
that \(Inv(N,f)=S\). The neighborhood \(N\) is called an {\em
isolating neighborhood of \(S\)}.

An {\em isolating block} \(N\) is a compactum such that
\(cl(int(N))=N\) and \(f^{-1}(N) \cap N \cap f(N) \subset
int(N)\). Isolating blocks are a special class of  isolating
neighborhoods.

 We consider the {\em exit set of} \(N\) to be defined as

\[
N^{-}=\{x \in N: f(x) \notin int(N)\}.
\]

Let \(S\) be an isolated invariant set and suppose \(L \subset N\)
is a compact pair contained in the interior of the domain of
\(f\). The pair \((N,L)\) is called a {\em filtration pair} for
\(S\) (see Franks and Richeson paper \cite{FR})  provided \(N\)
and \(L\) are each the closure of their interiors and

1) \(cl(N \setminus L)\) is an isolating neighborhood of \(S\),

2) \(L\) is a neighborhood of \(N^{-}\) in \(N\) and

3) \(f(L) \cap cl(N \setminus L) = \emptyset\).

\medskip

\begin{Remark}  Filtration pairs are easy to construct once we have an
isolating block \(N\). In fact, for every small enough closed
neighborhood \(L\) of \( N^-\), \((N,L)\) is a filtration pair
(\cite{FR}).
\end{Remark}

\medskip

In \cite{RS2} we compute the indices of the iterates of \({\Bbb
R}^3\)-homeomorphisms, \(f\), when there is a block \(N\), that is
topological closed ball, such that \(Inv(N,f)=\{0\}\). On the
other hand, there are not techniques for the explicit computation
of the sequence of the iterates of arbitrary homeomorphisms. Since
in this paper we shall deal with homeomorphisms such that for
every closed ball, \(B\), centered in \(0\), \(Inv(B,f) \cap
\partial B \ne \emptyset\), we will compute the sequence by
approximating adequately our map by a sequence of homeomorphisms
\(\{f_m\}_{m \in {\Bbb N}}\) such that \(Inv(B,f_m)=\{0\}\) for
every \(B\) and every \(m\in {\Bbb N}\).

In the following examples there will be an isolating block \(N\),
which is a solid ball, such that  \(Inv(N,f)=\{0\}\).  The
sequences of indices are easily seen to be periodic. However, they
will provide some ingredients we shall need to prove Theorems 1
and 2.

\medskip

\noindent {\bf 2.1 Examples where \(Inv(N,f)=\{0\}\)}.

\medskip

{\bf 1.} Consider  the linear homeomorphism \(g:{\Bbb R}^3
\rightarrow {\Bbb R}^3\) given by the matrix

\[
\begin{pmatrix}
1/2& 0& 0  \\
0&1/2& 0  \\
0 & 0 & 2
\end{pmatrix}.
\]

In this case, \(Fix(g)=\{0\}\), \(\{0\}\) is the unique compact
\(g\)-invariant set and \(0\) is an hyperbolic fixed point.

The computation of the sequence \(\{i(g^k,0)\}_{k \in {\Bbb N}}\)
is a very easy problem using standard methods. However, we are
going to calculate the sequence using different ideas.

\medskip

There exist a filtration pair \((N,L)\) such that \(N \) is an
isolating block, a closed 3-dimensional  ball, and \(L\) is a
disjoint union of two balls \(L_1\) and \(L_2\). Identifying
\(L_1\) and \(L_2\) to two different points \(q_1\) and \(q_2\) we
obtain a quotient space, denoted by \(N_L\), and a map induced by
\(g\), \(\bar{g}:N_L \rightarrow N_L\), with
\(Fix(\bar{g})=Per(\bar{g})=\{0, q_1, q_2\}\).

\medskip

Now, the Lefschetz number \[\Lambda(\bar{g}^k) =1=\]

\[=i_{N_L}(\bar{g}^k,0)+ i_{N_L}(\bar{g}^k,q_1) +
i_{N_L}(\bar{g}^k,q_2).\]

\medskip

Since \(q_1\) and \(q_2\) are attractors in \(N_L\), we have that
\(i_{N_L}(\bar{g}^k,q_1)=i_{N_L}(\bar{g}^k,q_2) =1\). Then
\(i(g^k,0)= i_{N_L}(\bar{g}^k, 0)= 1-2=-1\) for every \(k \in
{\Bbb N}\).

\medskip

On the other hand, since \(g\) preserves orientation,
\(i(g^{-k},0)= - i(g^k,0)=1 \) for every \(k \in {\Bbb N}\).

\medskip

Using similar ideas we can compute the sequences of indices of the
iterations of  \(g^{-1}\) in another way.  Analogously there
exists a filtration pair \((N, E)\) for \(g^{-1}\) such that \(N\)
is again an isolating block, a closed 3-dimensional ball and \(E\)
is now a solid torus.

\medskip

Identifying \(E\) to a point \(q\) we obtain the quotient space
\(N_E\) which is an ANR having the homotopy type of a 2-sphere.

\medskip

Now, \(2= \Lambda((\overline{g^{-1}})^{k}) =
i_{N_E}((\overline{g^{-1}})^{k},0)+
i_{N_E}((\overline{g^{-1}})^{k},q) \).

\medskip

Then, \(i_{{\Bbb R}^3}(g^{-k},0)= i_{N_E}((\overline{g^{-1}})^{k},
0)= 2-1=1\) for every \(k \in {\Bbb N}\), and again \(i_{{\Bbb
R}^3}(g^{k},0)= - i_{{\Bbb R}^3}(g^{-k},0)=-1 \) for every \(k \in
{\Bbb N}\).

\medskip
\medskip

 {\bf 2.} Let \(C\) be the cube \(C=[-1,1]^3\). Joining \(0\) with the one
skeleton of \(\partial C\) we obtain six closed and bounded cones
with disjoint interiors \(\{ c_j: j=1, \dots, 6\} \). Let \(C_j
=\{\lambda p: \lambda \geq 0 \text{ and } p \in c_j\}\). It is
clear that \({\Bbb R}^3= \bigcup_{j \in \{1, \dots, 6\}} C_j\).

\medskip

Let \(g\) be the homeomorphism of the first example. Let \(\pi^+=
\{(x,y,z) \in {\Bbb R}^3 : z \geq 0\}\). The restriction
\(g_{|\pi^+}\) is conjugated to a homeomorphism \(g_j: C_j
\rightarrow C_j\). Consider the orientation preserving
homeomorphism \( \phi: {\Bbb R}^3 \rightarrow {\Bbb R}^3\) defined
as \(\phi_{|C_j}= g_j\) for \(j \in \{1, \dots, 6\}\).

\medskip

We have that \(Per(\phi) = Fix(\phi)=\{0\}\) and, again, \(\{0\}\)
is the unique compact \(\phi\)-invariant set. Now, the stable
"manifold" is the cone of the one-dimensional skeleton of
\(\partial C\) and the unstable manifold decomposes into six one
dimensional branches (the union of the half lines joining \(0\)
with the center of each face of \(\partial C\)).

\medskip

It is not difficult to check that there exists a filtration pair
\((N,L)\) such that \(N\) is an isolating block 3-ball and \(L\)
is disjoint union of six 3-balls.

\medskip

If we identify each component of \(L\) to a different point we
have  the space \(N_L\) and the induced map \(\bar{\phi}: N_L
\rightarrow N_L\). Now, \(Per(\bar{\phi})=Fix(\bar{\phi})= \{0,
q_1, q_2, \dots, q_6\}\). All \(q_j\)'s are attractors and then
they have index =1.

\medskip

Therefore, \(i(\phi^{k},0)= 1-6=-5\) for every \(k \in {\Bbb N}\).
As a consequence, \(i(\phi^{-k},0)= - i(\phi^k, 0)= 5\) for every
\(k \in {\Bbb N}\).

\medskip

We can  compute the sequence of indices of \(\phi^{-1}\) directly
by considering a filtration pair for \(\phi^{-1}\). Indeed, there
is a pair \((N,E)\) where \(N\) is a closed ball and \(E\) is an
adequate tubular neighborhood of \(\partial N \cap (\bigcup_{j \in
\{1, \dots, 6\}} \partial C_j) \). Now, the quotient space \(N_E\)
is an ANR having the homotopy type of the wedge of five 2-spheres.
Each  of these five 2-spheres corresponds to one of the faces of
\(\partial C\). The remaining one represents the sum of the
others.

\medskip

Now, it is easy to check, by choosing obvious generators, that the
matrix of \(\overline{\phi^{-1}}_{2}^{*}: H_2(N_E) \rightarrow
H_2(N_E)\) can be assumed to be the identity.

\medskip

Then, the Lefschetz number \(\Lambda((\overline{\phi^{-1}})^{k})=
6\) for all \(k \in {\Bbb N}\) and \(i(\phi^{-k},0)=6-1= 5\), \(k
\in {\Bbb N}\).

\medskip

\medskip

{\bf 3.} Let \(\psi: {\Bbb R}^3 \rightarrow {\Bbb R}^3\) be the
homeomorphism given by the composition of \(r_{\pi/2} \circ
\phi^{-1}\) where \(\phi \) is the homeomorphism of Example 2 and
\(r_{\pi/2}\) is the \(\pi/2\)-rotation with respect to the axis
\(\{(x,y,z) \in {\Bbb R}^3 : x=y=0 \}\).

Here, the sequence is periodic of period 4 and it is not difficult
to show that

\[ i(\psi^{-k}, p)=\left\{\begin{array}{ll}
-5 & \text{ if } k \in 4{\Bbb N} \\
-1 & \text{ if } k \notin 4{\Bbb N} \end{array} \right.
\]

\noindent and

\[ i(\psi^{k}, p)=\left\{\begin{array}{ll}
5 & \text{ if } k \in 4{\Bbb N} \\
1 & \text{ if } k \notin 4{\Bbb N} \end{array} \right.
\]

\medskip
\medskip
\medskip
\medskip

\centerline{{\bf 3. The construction of the homeomorphisms. Proof
of the Theorems. }}

\medskip
\medskip

\noindent {\bf Proof of Theorem 1.}

\medskip
There is no lost of generality  if we assume that \(\{c_m\}_{m \in
{\Bbb N}} \subset {\Bbb N}\).

Since we will be interested just in the elements of the sequence
\(\{c_m\}_{m \in {\Bbb N}}\) with \(m\) prime, we shall rename
some of the terms of that sequence in the following way. If \(q
\in {\Bbb N}\) is the \(k\)-th prime number, we will write
\(c_q=c'_k\).

Our aim is to construct a ${\mathbb R}^3$-homeomorphism $f$ such
that $Fix(f)=Per(f)=\{0\}$ and $\limsup
\frac{|i(f^m,0)|}{c_m}=\infty$. For this end it is enough that for
each  \(k\), if \(p\) denotes the \(k\)-th prime, we get the index
of $f^p$ to be

\[ i(f^p,0)=-p^{c_p}=-p^{c'_k}\].

To simplify the notation we will write again the new sequence
\(\{c'_{m}\}_{m \in {\Bbb N}} \) as \(\{c_{m}\}_{m \in {\Bbb
N}}\).

Let $N=B(0, 1)$ be the unit closed ball. Our first step is to make
a partition of $N$ in solid regions. Each of these regions will
have a different dynamics.

Let $A_i \subset N$ be the solid region limited by a cone $C_i$
with vertex $0$ and axis the line joining the poles $n$ and $s$ of
$N$ in such a way that

\[
A_i \subsetneq A_{i+1} \subsetneq \cdots \subset N^+ \quad \text{ and }
\quad cl(\bigcup A_i)=N^+
\]

\noindent with $N^+=\{\bar x \in N : x_3 \geq 0\}$.

We define the different solid regions on which we will have the
characteristic dynamics of $f$ in the next way:

Let $S_0=A_0$, $S_i=cl(A_i \setminus A_{i-1})$ for $i \in {\mathbb N}$ and
let $S_{\infty}=\{\bar x \in N : x_3 \leq 0\}=N^-$.

We have a decomposition of $N$,

\[
N=\bigcup_{m=0}^{\infty} S_m
\]

\noindent with $S_i \cap S_{i+1}=C_i \cap N^+$ and $S_i \cap
S_j=\{0\}$ if $j \notin \{i-1,i,i+1\}$.

We construct the sets $S_i$, $i \notin \{0, \infty\}$, in such a
way that the length of each of the two arcs $c_i \cup
c'_i=\partial(N) \cap S_i \cap \{x=0\}$ is
$l_i=\frac{\pi}{2^{i+2}}$. The length of the arc $c_0=\partial(N)
\cap S_0 \cap \{x=0\}$ is $\frac{\pi}{2}$. Let us observe that if
we work in spherical coordinates $(\rho, \theta,\phi)$, the angle
$\phi$ of the points in $C_n^+=C_n \cap \pi^+$ is

\[
\phi=\sum_{i=0}^{n}
\frac{\pi}{2^{i+2}}=\frac{\pi}{2}-\frac{\pi}{2^{n+2}}
\]

We will define the homeomorphism $f$ as the composition of two
homeomorphisms $f_0:{\Bbb R}^3 \rightarrow {\Bbb R}^3$ and
$g_0:{\Bbb R}^3 \rightarrow {\Bbb R}^3$ with $f=f_0 \circ
g_0:{\Bbb R}^3 \rightarrow {\Bbb R}^3$.

Let

\[
E(A_i)=\{\lambda \bar x \, : \, \lambda \in {\mathbb R}^+ , \bar x
\in A_i\}
\]

In the same way we define the sets $E(S_i)$ with $i=0, \dots,
\infty$.

Let $\{r_m\}_{m \in {\Bbb N}}=\{p_m/q_m\}_{m \in {\Bbb N}}$ be a
sequence of rational numbers converging to an irrational number
$r$ with $\{q_m\}_{m \in {\Bbb N}}$ the sequence of prime numbers
and such that $g.c.d.(p_m,q_m)=1$. We can construct the sequence
$\{r_m\}_{m \in {\Bbb N}}$ with $0<r<1$ in the following way:

For each $q_m$ we consider a partition of the unit interval
$[0,1]$ in $q_m$ intervals of length $1/q_m$ and select $p_m <
q_m$ as the natural number such that $d(p_m/q_m, r)=\min
\{d(n/q_m, r)\}$ with $n \in {\mathbb N}$. Then, the sequence
$\{p_m/q_m\}_{m \in {\Bbb N}} \to r$ when $m \to \infty$ and
$g.c.d(p_m, q_m)=1$.


In $S_n$ with $n=2m-1$ odd we consider a family of $q_m^{c_m}$
isometric solid regions $\{T_{j,m}\}$, linearly isomorphic to the
sets $A_i$, and such that $T_{j,m} \subset S_n$, $T_{j,m} \cap
\partial(S_n)=D_{j,m} \cup \{\bar 0\}$ with $D_{j,m}$ a closed disc.
We put these solid regions in $S_n$ with constant angle
$\frac{2\pi}{q_m^{c_m}}$ around the vertical axis (which joins the
poles of $N$) and define $E(T_{j,m})=\{\lambda \bar x \, : \,
\lambda \in {\mathbb R}^+, \bar x \in T_{j,m}\}$. See figure 1.

\medskip
\medskip

\centerline{
\includegraphics[width=9cm,height=9cm,keepaspectratio,clip]{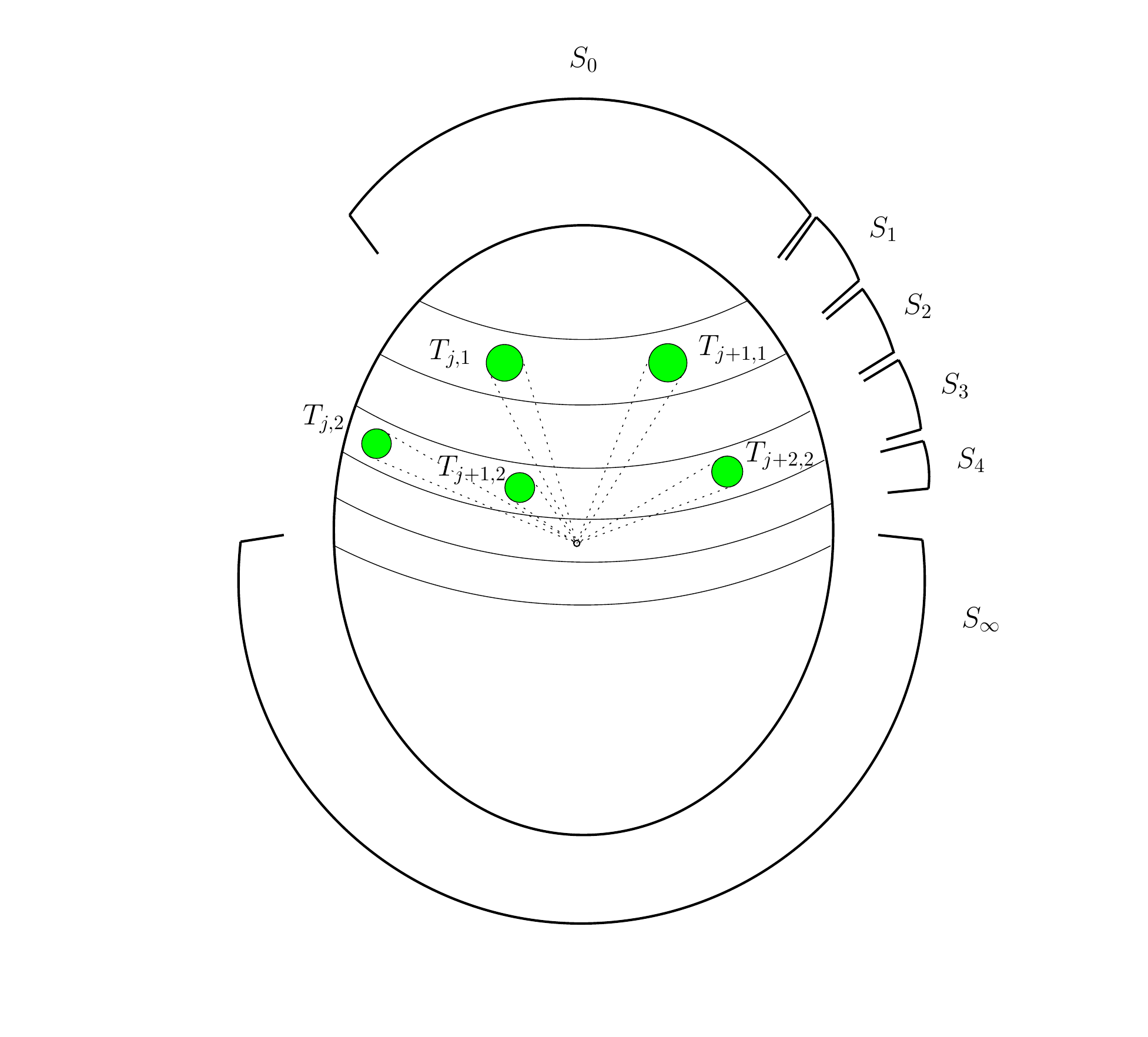}
}

\centerline{Figure 1}

\medskip
\medskip


The dynamics of $f_0$ in $E(S_{\infty})$ is

\[
f_0(\bar x)=(1-\frac{1}{\pi}\lambda(\bar x)) \bar x
\]

\noindent with $\lambda(\bar x)$ defined as the length of the
parallel arc of $\partial(N)$ joining $\frac{\bar x}{|\bar x|}$
with the plane $\{z=0\}$. Working with spherical coordinates we
have

\[
f_0(\rho, \theta, \phi)=\left( \left(
\frac{1}{2}+\frac{-\phi}{\pi} \right) \rho, \theta, \phi \right)
\]


The dynamics of $f_0$ in the two consecutive conical regions
$C_{2m}^+ \cup C_{2m+1}^+$ is

\[
f_0(\rho, \theta,\phi)=\left( \left( 1-\frac{1}{2^{2m}}\right)
\rho, \theta, \phi \right)
\]


In the same way, the dynamics of $f_0$ in each region
$cl(E(S_{2m+1}) \setminus \bigcup E(T_{j,m+1}))$ is

\[
f_0(\rho, \theta, \phi)=\left( \left( 1-\frac{1}{2^{2m}} \right)
\rho, \theta,\phi \right)
\]

\medskip
\medskip

\vspace{0.3cm}

\centerline{
\includegraphics[width=7cm,height=7cm,keepaspectratio,clip]{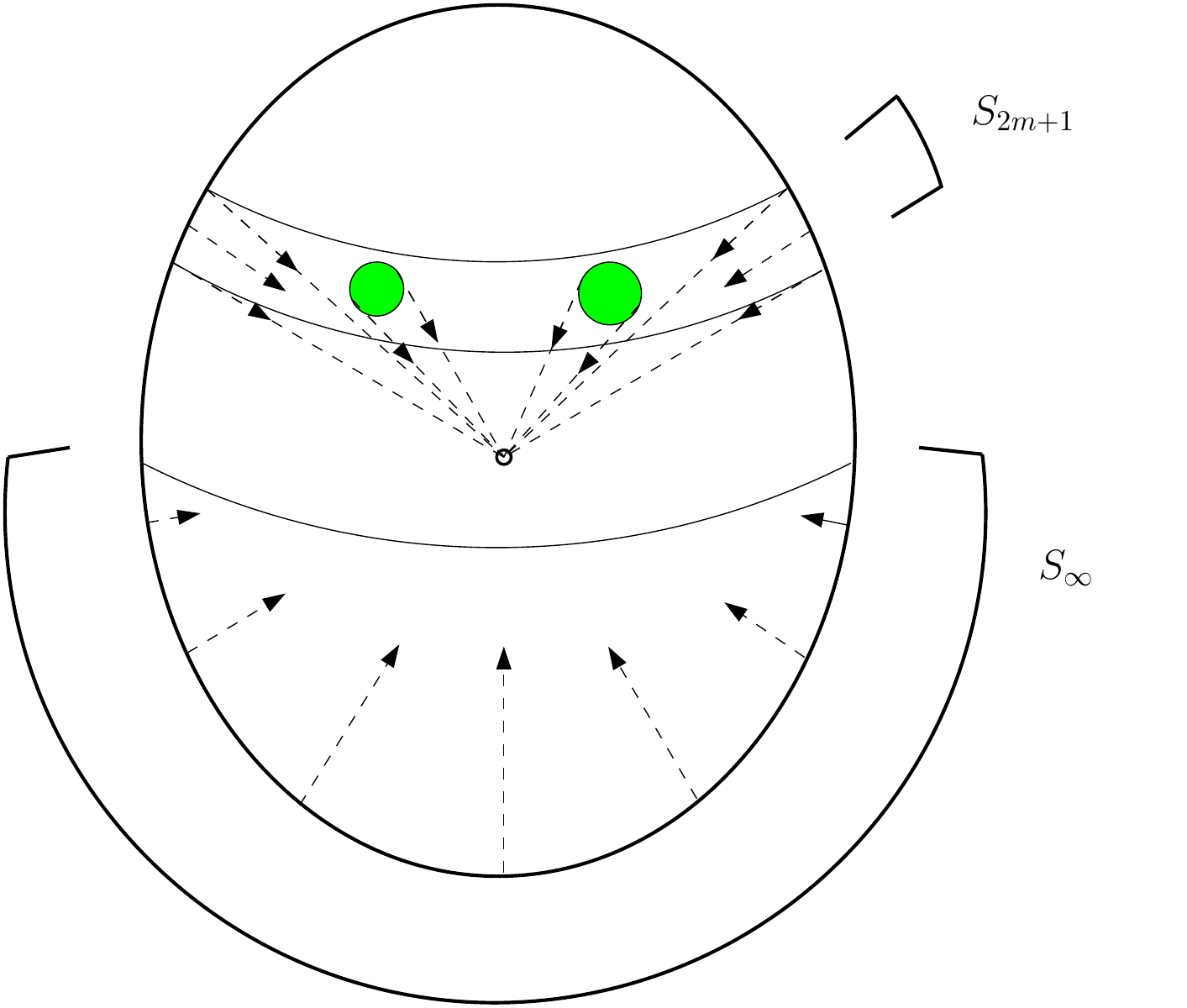}
}

\centerline{Figure 2}

\medskip
\medskip

On the other hand, the dynamics of $f_0$ in the sets $E(T_{j,m})
\subset S_n$ for $n=2m-1$ is conjugated with the given in Example
1 for the map $g|_{\pi^+}:\pi^+ \rightarrow \pi^+$ and commutes
with a rotation of angle $\frac{2\pi}{q_m^{c_m}}$. Moreover, we
construct $f_0$ in $E(T_{j,m}) \subset S_n$ in such a way that
$d(f_0(\bar x), \bar x) \leq k_n \|\bar x\|$ for all $\bar x \in
E(T_{j,m})$ and with $k_n \to 0$ when $n \to \infty$. See figure
3.

\medskip
\medskip

\centerline{
\includegraphics[width=5cm,height=5cm,keepaspectratio,clip]{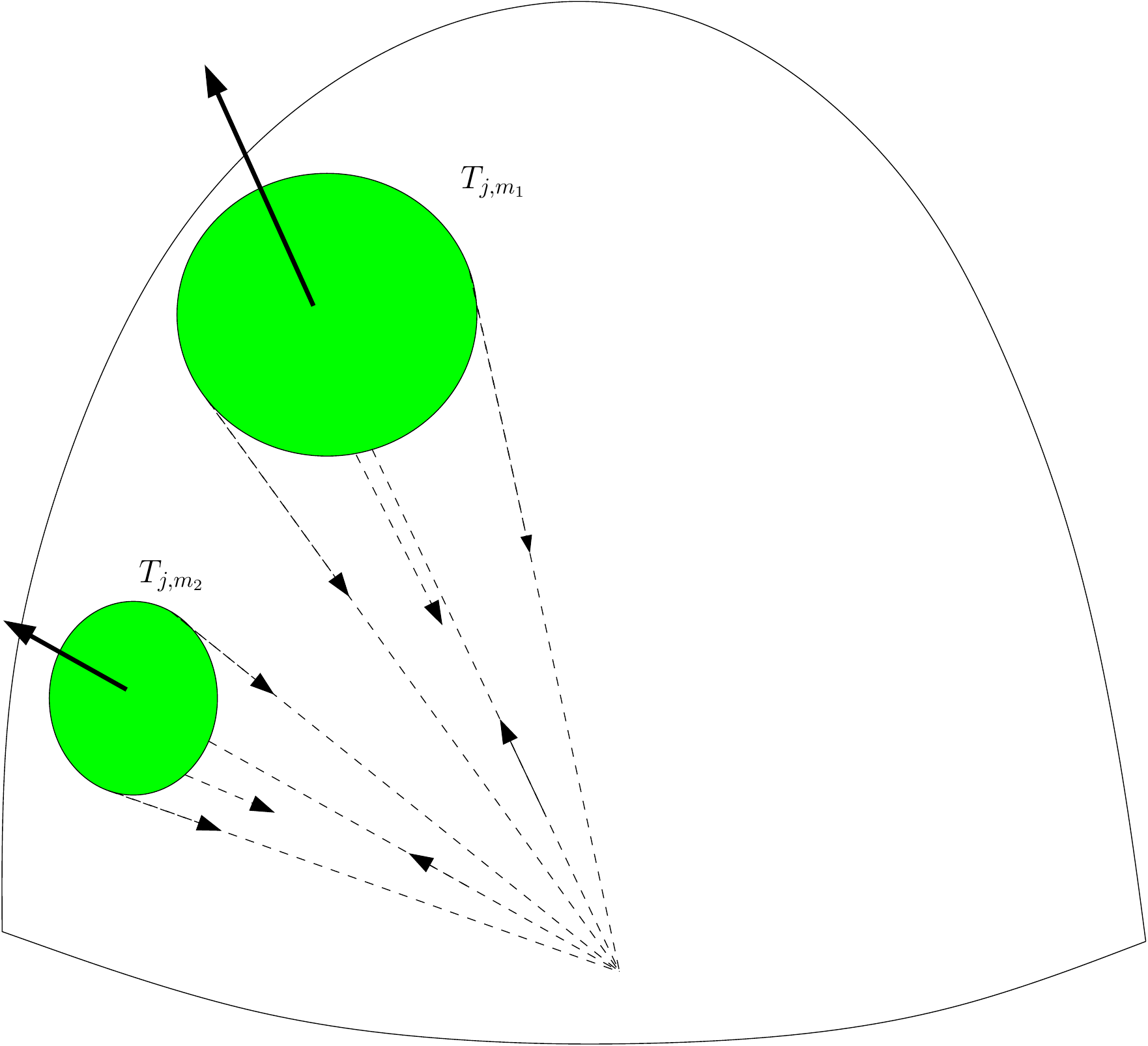}
}

\centerline{Figure 3}

\medskip
\medskip


Let us suppose that $n=2m$ even. For every point $\bar x \in
E(S_{2m})$, the coordinate $\phi$ is in the interval $\left[
\frac{\pi}{2}-\frac{\pi}{2^{2m+1}},
\frac{\pi}{2}-\frac{\pi}{2^{2m+2}} \right]$. The dynamics of $f_0$
in the regions $E(S_{n})$ with $n=2m$ is, taking spherical
coordinates,

\[
f_0(\rho, \theta,\phi)=(k_n(\phi) \rho, \theta, \phi)
\]

\noindent with

\[
k_n:\left[
\frac{\pi}{2}-\frac{\pi}{2^{2m+1}},\frac{\pi}{2}-\frac{\pi}{2^{2m+2}}
\right] \rightarrow \left[ 1-\frac{1}{2^{2(m-1)}},
1-\frac{1}{2^{2m}} \right]
\]

\noindent an increasing, bijective linear map.

\medskip
\medskip

\centerline{
\includegraphics[width=7cm,height=7cm,keepaspectratio,clip]{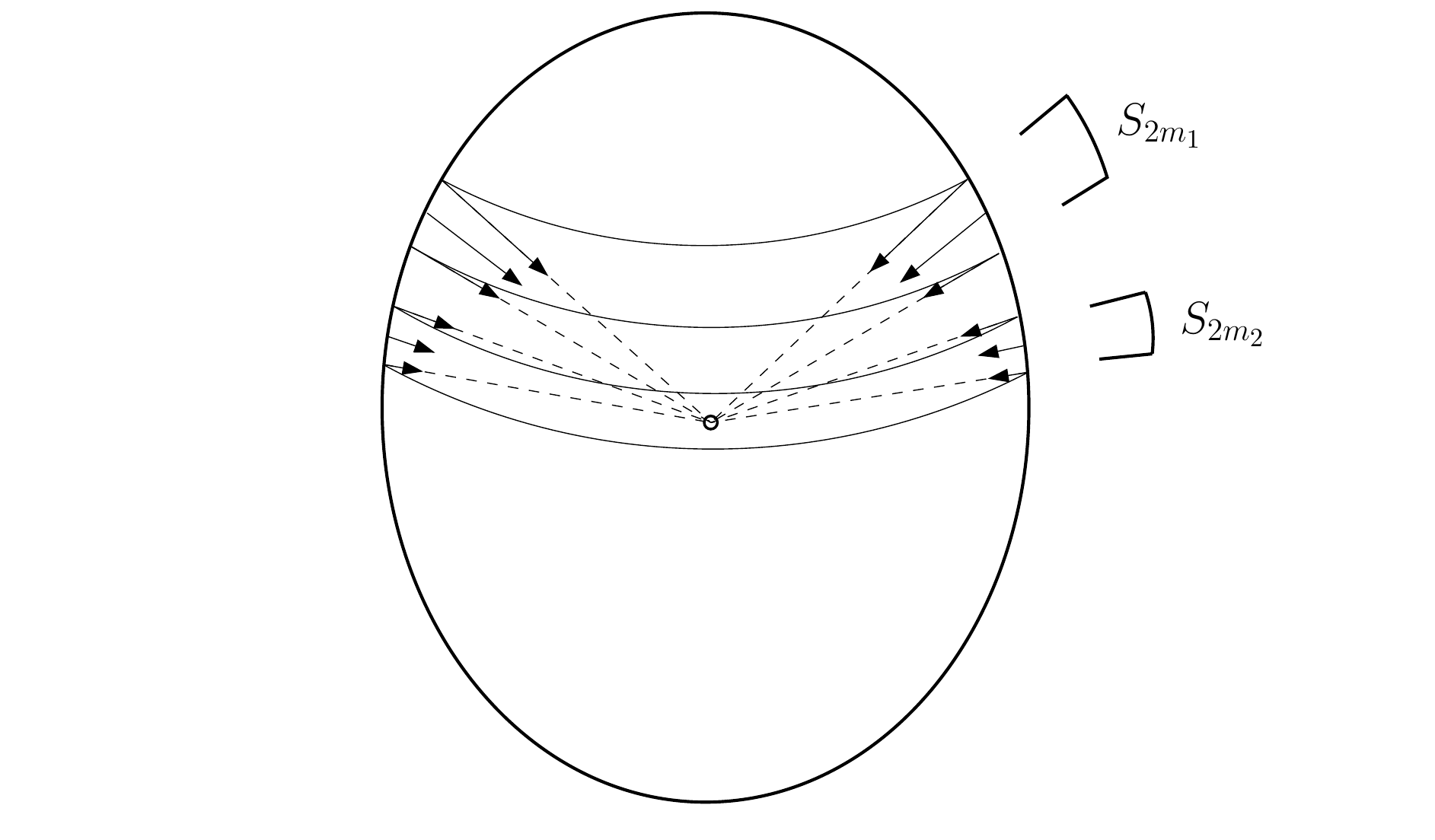}
}

\centerline{Figure 4}

\medskip
\medskip


For each solid region $T_{j,m}$, the exit set of $f_0|_{T_{j,m}}$
is a closed ball $L_{j,m}$ such that $L_{j,m} \cap \partial(N)$ is
a closed disc. These closed balls are the exit regions of $N$ for
$f_0|_{S_n}$ and have constant angle $\frac{2\pi}{q_m^{c_m}}$
around the vertical axis (which joins the poles of $N$). See
figure 5

\medskip
\medskip

\centerline{
\includegraphics[width=8cm,height=8cm,keepaspectratio,clip]{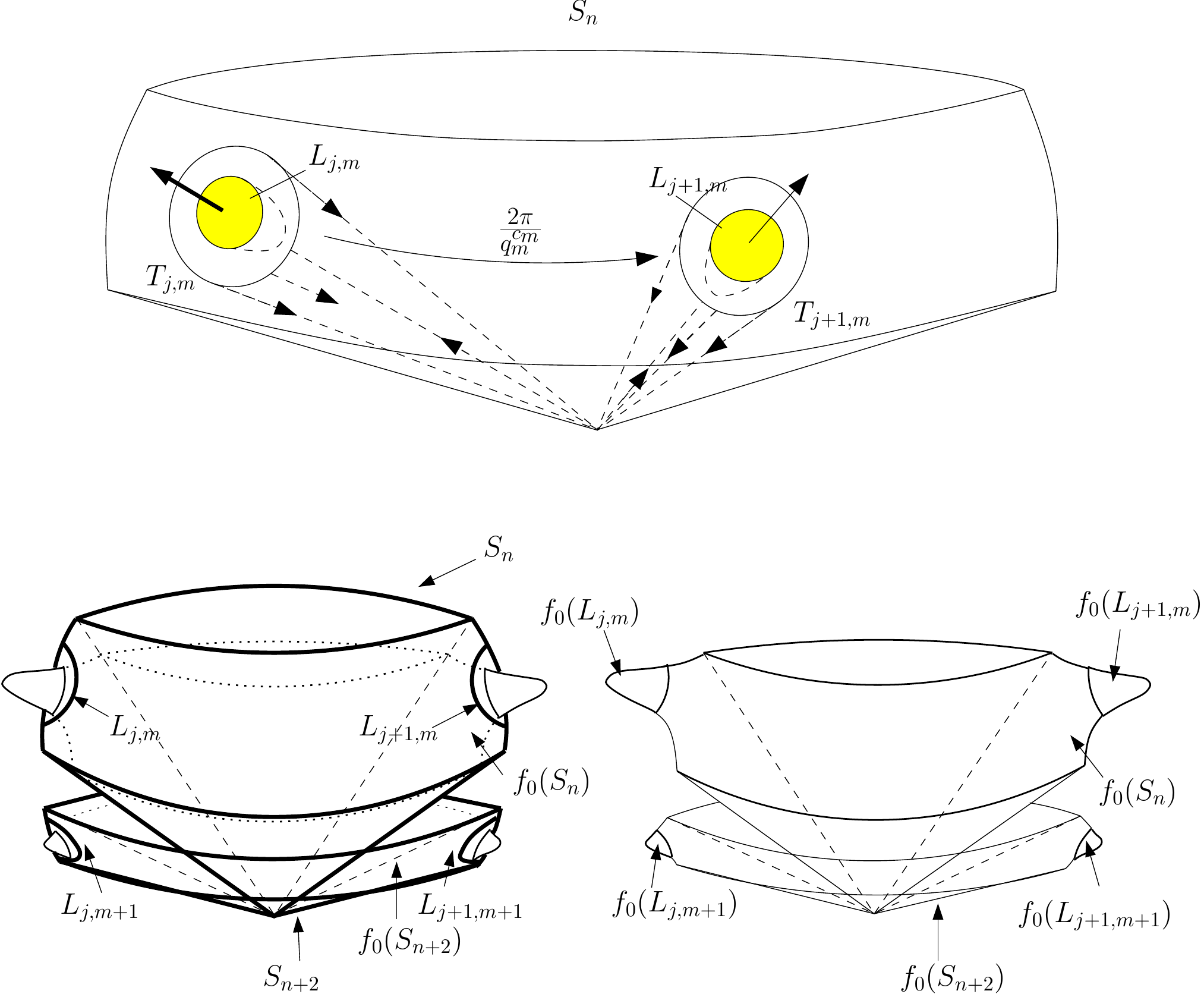}
}

\centerline{Figure 5}

\medskip
\medskip

It only remains to construct the dynamics of $f_0|_{E(S_0)}$. It
is topologically conjugated with the given in Example 1 for
$g|_{\pi^+}$. We obtain an exit region for $f_0|_{S_0}$ which is a
closed ball $L_0 \subset S_0$ such that $L_0 \cap
\partial(N)$ is a closed disc. The dynamical behavior is equivalent
to the dynamics obtained in the sets $T_{j,m}$. See figure 6.

\medskip
\medskip

\centerline{
\includegraphics[width=9cm,height=9cm,keepaspectratio,clip]{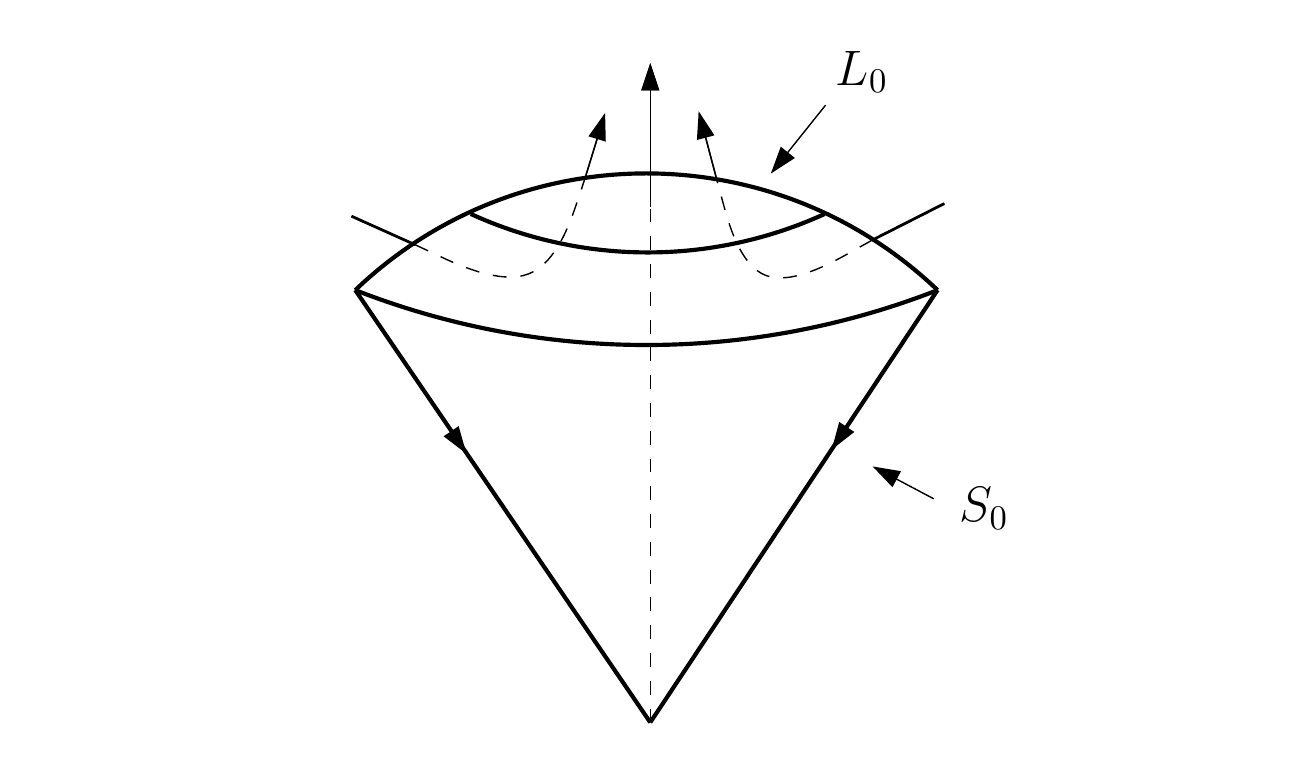}
}

\centerline{Figure 6}

\medskip
\medskip

It is easy to check that the map $f_0:{\Bbb R}^3 \rightarrow {\Bbb
R}^3$ is a homeomorphism, limit of homeomorphisms $\{f_{0,n}\}_n$,
with $f_{0,n}:{\Bbb R}^3 \rightarrow {\Bbb R}^3$ defined as

\[
f_{0,n}(\bar x)=\left\{\begin{array}{ll} f_0(\bar x) & \text{ if }
\bar x \in E(A_n) \cup E(A_n^-) \\
k(n) \bar x & \text{ if } \bar x \in N \setminus (E(A_n) \cup
E(A_n^-))
\end{array}\right.
\]

\noindent where $A_n^-=\{\bar x \in N \text{ such that } -\bar x
\in A_n\}$ and $k(n)=1-\frac{1}{2^{n-1}}$ if $n$ is odd and
$k(n)=1-\frac{1}{2^{n}}$ if $n$ is even.

Let us observe that
$Fix(f_{0,n})=Per(f_{0,n})=Inv(N,f_{0,n})=\{0\}$ and
 $Fix(f_0|_N)=Per(f_0|_N)=Inv(N,f_0)=N \cap \{z=0\}$ with $N^-(f_0)=\bigcup L_{j,m} \cup
(\{z=0\} \cap \partial(N)) \cup L_0$.


The homeomorphism $g_0:{\Bbb R}^3 \rightarrow {\Bbb R}^3$ is
defined in the next way:

The map  $g_0|_{E(S_{n})}$ with $n=2m-1$ odd is a rotation around
the vertical axis with angle $2\pi \frac{p_m}{q_m}$, that is,

\[
g_0|_{E(S_{2m-1})}(\rho, \theta\, \phi)=(\rho, \theta+2\pi
\frac{p_m}{q_m}, \phi)
\]

The restrictions $g_0|_{E(S_{\infty})}$ and $g_0|_{E(S_0)}$ are
rotations around the vertical axis with angles $2\pi r$ and $2\pi
\frac{p_1}{q_1}$ respectively.

The dynamics of $g_0|_{E(S_n)}$ with $n=2m$ even is the following:

Since $g_0|_{C_{n-1} \cap \pi^+}$ and $g_0|_{C_{n} \cap \pi^+}$
are rotations with angles $2\pi \frac{p_m}{q_m}$ and $2\pi
\frac{p_{m+1}}{q_{m+1}}$, given a cone $C$ with vertex $0$ and
axis the line joining the poles of $N$ such that $C \cap \pi^+
\subset E(S_n)$, we construct the dynamics in $C \cap \pi^+$ as a
rotation with angle $c \in \left[ 2\pi \frac{p_m}{q_m}, 2\pi
\frac{p_{m+1}}{q_{m+1}} \right]$ in such a way that $c$ tends to
$2\pi \frac{p_m}{q_m}$ ($2\pi \frac{p_{m+1}}{q_{m+1}}$) if $C$
tends to $C_{n-1}$ ($C_{n}$). Working with spherical coordinates

\[
g_0|_{S_{2m}}(\rho, \theta\, \phi)=(\rho, \theta+k_n(\phi), \phi)
\]

\noindent with

\[
k_n:\left[ \frac{\pi}{2}-\frac{\pi}{2^{2m+1}},
\frac{\pi}{2}-\frac{\pi}{2^{2m+2}}\right] \rightarrow \left[
2\pi\frac{p_m}{q_m}, 2\pi\frac{p_{m+1}}{q_{m+1}} \right]
\]

\noindent an increasing, bijective linear map.

\medskip
\medskip

\centerline{
\includegraphics[width=9cm,height=9cm,keepaspectratio,clip]{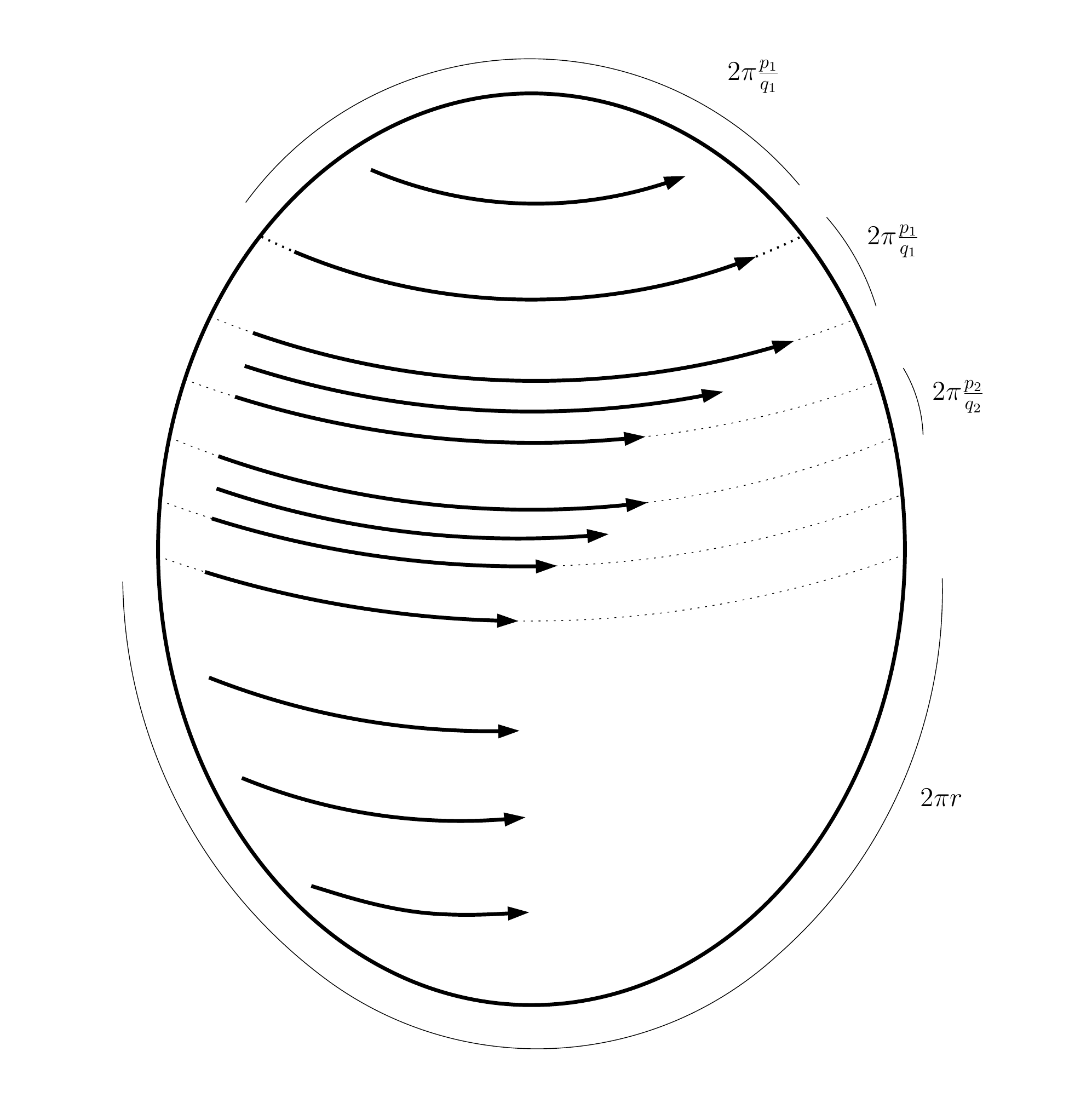}
}

\centerline{Figure 7}

\medskip
\medskip

The map $g_0:{\Bbb R}^3 \rightarrow {\Bbb R}^3$ constructed is a
homeomorphism and limit of homeomorphisms $\{g_{0,n}\}_n$ with
$g_{0,n}$ defined as follows:

Let us define $g_{0,n}$ for $n=2m-1$ odd (if $n$ is even, the
construction is analogous). Given $\bar x \in E(A_n) \cup
E(A_n^-)$ we define $g_{0,n}(\bar x)=g_0(\bar x)$. On the other
hand, let us observe that for every $\bar x \in C_n \cap \pi^+$
the spherical coordinates are $(\rho, \theta, \phi_n)$ with
$\phi_n=\frac{\pi}{2}-\frac{\pi}{2^{n+1}}$ fixed. Since
$g_0|_{C_{n} \cap \pi^+}$ and $g_0|_{C_n \cap \pi^-}$ are
rotations with angles $2\pi \frac{p_{m}}{q_{m}}$ and $2\pi r$, for
each $\bar x \in cl({\mathbb R}^3 \setminus (E(A_n) \cup
E(A_n^-)))$ with spherical coordinates $(\rho, \theta , \phi)$ we
construct the dynamics of $g_{0,n}$ as a rotation with angle
$c(\phi) \in \left[ 2\pi \frac{p_{m}}{q_{m}}, 2\pi r \right]$ in
such a way that $c$ tends to $2\pi \frac{p_{m}}{q_{m}}$ ($2\pi r$)
if $\phi$ tends to $\phi_n$ ($\frac{\pi}{2}-\phi_n$). Then, given
$n=2m-1$,

\[
g_{0,n}|_{cl({\mathbb R}^3 \setminus (E(A_n) \cup E(A_n^-)))}
(\rho, \theta, \phi)=(\rho, \theta+k_n(\phi),\phi)
\]

\noindent with

\[
k_n: \left[ \frac{\pi}{2}-\frac{\pi}{2^{2m}},
\frac{\pi}{2}+\frac{\pi}{2^{2m}} \right] \rightarrow \left[2\pi
\frac{p_m}{q_m}, 2\pi r \right]
\]

\noindent an increasing, bijective linear map.

The map $f=f_0 \circ g_0: {\Bbb R}^3 \rightarrow {\Bbb R}^3$ is a
homeomorphism with $Fix(f)=Per(f)=\{0\}$ and $Inv(N,f)=N \cap
\{z=0\}$. If we consider the sequence of homeomorphisms
$\{f_n\}_n$ with $f_n=f_{0,n} \circ g_{0,n} :{\Bbb R}^3
\rightarrow {\Bbb R}^3$ we have

\[
Fix(f_n)=Per(f_n)=Inv(N,f_n)=\{0\}
\]

\noindent and it is obvious that $f$ is limit of the
homeomorphisms $\{f_n\}$.

Let us compute the fixed point index $i(f^n,0)$ for $n \in
{\mathbb N}$. For this purpose we will use the next two results of
existence of homotopies between near enough maps and homotopy
invariance of the fixed point index.

\begin{Remark}
Let $f:X \subset {\mathbb R}^m \rightarrow {\mathbb R}^n$ be a
continuous map. Then if $g:X \rightarrow {\mathbb R}^n$ is a
continuous map near enough $f$, they are homotopic.
\end{Remark}

\begin{Remark} Let $X$ be a metric ANR, $W$ an open subset of $X$ and $F:cl(W) \times [0,1] \rightarrow X$ a continuous and compact map such
that $F(x,t) \neq x$ for $(x,t) \in \partial(W) \times [0,1]$.
Then $i_X(F_t, W)$ is constant for $0 \leq t \leq 1$.
\end{Remark}

Let us fix $d \in {\mathbb N}$. Since the map $f^d|_N: N
\rightarrow {\mathbb R}^3$ can be approximated by maps of the type
${f_n}^d|_N: N \rightarrow {\mathbb R}^3$, from the first of the
two remarks there exists $n_0 \in {\mathbb N}$ such that for each
$n \geq n_0$ there exists a homotopy $H: N \times I \rightarrow
{\mathbb R}^3$ with $H_0=f^d$, $H_1={f_n}^d$ and $H(x,t) \neq x$
for all $x \in \partial(N)$ and $t \in [0,1]$. From the second
remark, we obtain that $i(f^d, 0)=i({f_n}^d, 0)$.

Let us compute $i({f_n}^d, 0)$. There exists a finite family of
closed balls $\{L_{j,m}\}$ contained in $N$ which are the exit
regions of $N$ for $f_n|_N$. Identifying the sets $\{L_{j,m}\}$ to
points $\{l_{j,m}\}$ we obtain a quotient space $N_L$, which is a
closed ball, and an induced map $\bar{f_n}: N_L \rightarrow N_L$.
It is obvious that $i_{N_L}(\bar{f_n}^d, 0)=i({f_n}^d,0)$. Given
$m$ fixed, the action of the map $\bar{f_n}$ on the family of
points $\{l_{j,m}\}_j$, with $j=1, \dots, q_m^{c_m}$, give us a
union of $q_m^{c_m-1}$ cycles of length $q_m$,

\[
\{l_{j,m}\}_j=\bigcup_k \{l(k,1), \dots, l(k,q_m)\},
\]

\noindent with $k=1, \cdots, q_m^{c_{m-1}}$, such that

\[
\bar{f_n}(l(k,r))=l(k,r+1)
\]

\noindent for $r=1, \dots, q_m$.

\medskip
\medskip

\centerline{
\includegraphics[width=9cm,height=9cm,keepaspectratio,clip]{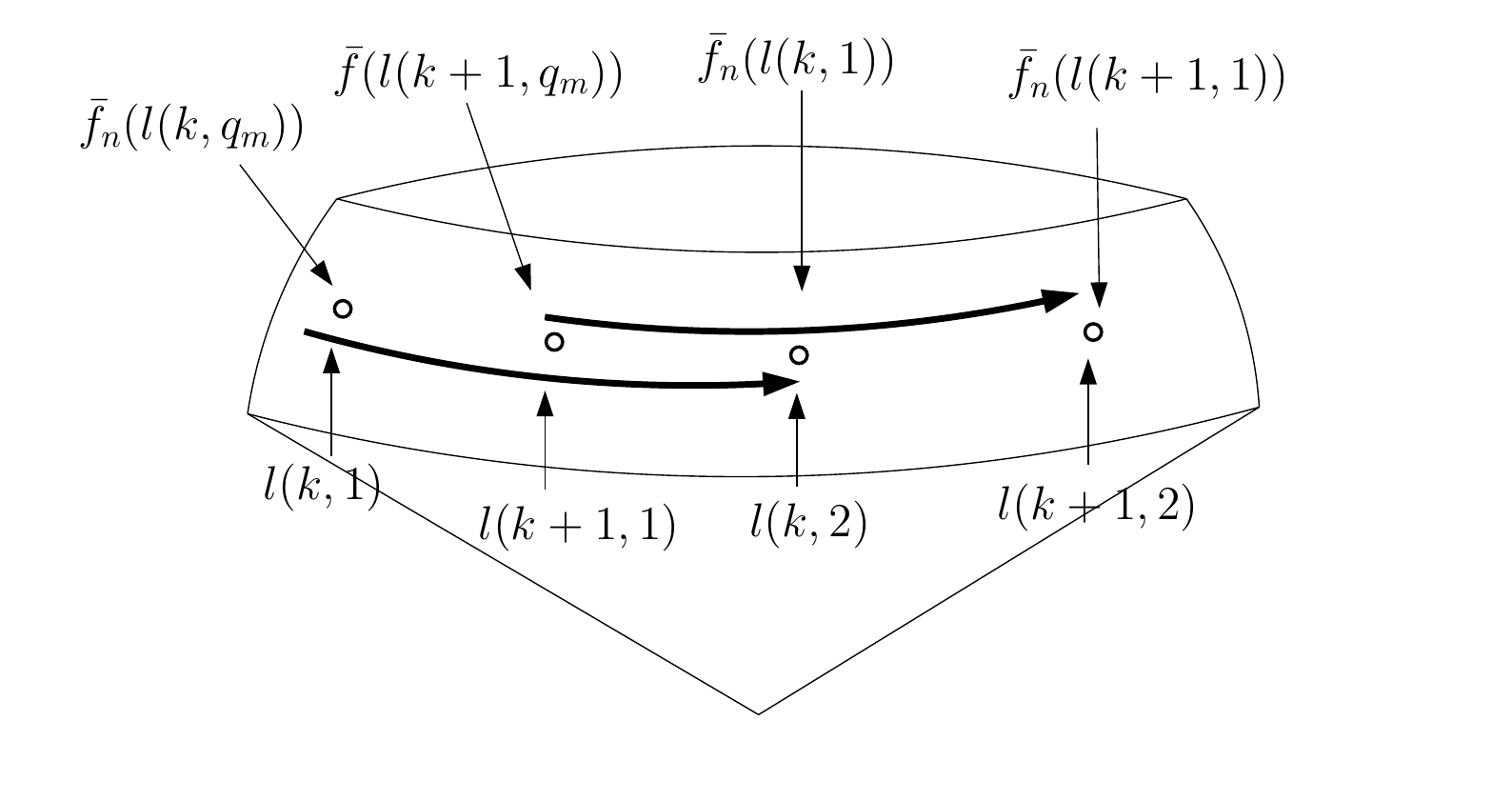}
}

\centerline{Figure 8}

\medskip
\medskip

It is obvious that

\[
i_{N_L}(\bar{f_n}^d, l_{j,m})=\left\{\begin{array}{ll}
1 & \text{ if } d \in q_m {\Bbb N} \\
0 & \text{ if } d \notin q_m {\Bbb N}
\end{array}\right.
\]

We obtain the equality

\[
1=i_{N_L}(\bar{f_n}^d, N)=i(f^d, 0)+\sum_{\stackrel{j=1, \dots,
{q_m}^{c_m}}{q_m | d}} i_{N_L}(\bar{f_n}^d, l_{j,m})+1
\]

\noindent where the  last 1 is due to the dynamics in the closed ball
$L_0$. Then,

\[
i(f^d, 0)=-\sum_{\stackrel{j=1, \dots, {q_m}^{c_m}}{q_m | d}}
i_{N_L}(\bar{f_n}^d, l_{j,m})=-\sum_{q_m | d} {q_m}^{c_m}
\]

If we consider $d=q_m$ prime,

\[
i(f^{q_m}, 0)=-{q_m}^{c_m}
\]

\noindent and the result is proved.


On the other hand, let us observe that if we consider the sequence
of natural numbers $\{q_m^k\}_k$ with $q_m$ the \(m\)-th prime
number and $k \in {\mathbb N}$, then

\[
i(f^{q_m^k}, 0)=-{q_m}^{c_m} \text{ for all } k \in {\mathbb N}.
\]

Let us observe also that for each $m=p_1^{r_1} \cdots p_k^{r_k}
\in {\mathbb N}$, with $p_1, \dots, p_k$ different prime numbers,
$i(f^m, 0)=i(f_n^m, 0)=-\sum_{j=1, \cdots, k} {p_j}^{c_j}$ for
every $n \geq \max \{p_1, \dots, p_k\}$. \(\square\).

\medskip
\medskip

\begin{Remark} One can consider the dual construction of Theorem 1, i.e.
the map \(f\) at \(\infty\) and the inverse homeomorphism
\(f^{-1}\). In the first case, for every closed ball, \(B\),
centered in \(\infty\), \(Inv(B,f) \cap
\partial B \ne \emptyset\) and in  the latter  the exit sets for
each of the analogous approaching homeomorphisms, \(h_m\), are
solid \(( 1+ \sum q_{k_m}^{c_{k_m}})\)-tori. Following similar
arguments (see also the examples in Section 2), one has that

\[
i(f^{-n}, 0)=i(f^n, \infty)=\sum_{q_m | n} {q_m}^{c_m}.
\]

Consequently, if we see the homeomorphism \(f\) as a
\(S^3\)-homeomorphism such that \(Fix(f)=Per(f)=\{0, \infty\}\),
it follows that  \(\limsup \frac{|i(f^m, 0)|}{c_m}= \limsup
\frac{|i(f^m, \infty)|}{c_m}= \infty\).
\end{Remark}

\medskip
\medskip

\noindent {\bf Proof of Theorem 2.}
\medskip

The ingredients of the proof of Theorem 2 are the homeomorphisms
given in Theorem 1 and the plug construction developed by Wilson
in \cite{W} (see also \cite{BV}). We shall maintain the notation
of Theorem 1.

Consider the solid cylinder \(B= \{(x,y,z) \in {\Bbb R}^3: x^2 +
y^2 \leq 1, z \in [a,b]\}\) and the  flow induced by the constant
vector field \(Y=(0,0,1)\). Denote respectively by \(\sigma(B)\),
\(\tau(B)\) and \(\beta(B)\) the lateral, top and bottom
boundaries of \(B\).

A {\em flow box} \((U,g)\) for a vector field \(X\) at a point
\(p\) consists of a neighborhood \(U\) of \(p\) and a
diffeomorphism \(g:B \rightarrow U\) such that:

i) \(X\) is transverse to \(g(\beta(B))\).

ii) There is a positive constant \(c\) such that
\(\phi(ct,g(x))=g(\psi(t,x))\) where \(\phi(t,\cdot)\) and
\(\psi(t, \cdot)\) denote the flows induced by \(X\) and \(Y\) on
\(B\) respectively. When it is clear from the context, we shall
omit the diffeomorphism \(g\).

Let \(U\) and \(V\) be two flow boxes with \(V\subset U\). Then
\(V\) is called {\em  a shrinkage} of \(U\) if \( \sigma(V)
\subset int(U)\), \(\tau(V) \subset \tau(U)\) and \(\beta(V)
\subset \beta(U)\).

Let us recall the following version of Wilson's theorem (\cite{W})
that we will need.

\begin{Thm}
Let \(X\) be a  \(C^{\infty}\) \({\Bbb R}^3\)-vector field. Let
\(U\) be a flow box of \(X\) and let \(V\) be a shrinkage of
\(U\). Then, there exist a \(C^{\infty}\) vector field \(X^1\) on
\(U\) such that:

a) \(X^1\) coincides with \(X\) on a neighborhood of \(\partial
U\).

b) The limit sets of \(X^1\) are a finite collection of invariant
circles on which the restricted flow is minimal.

c) Every trajectory of \(X^1\) which intersects \(\beta(V)\)
remains in positive time inside \(U\).

d) Each trajectory of \(X^1\) which leaves \(U\) in positive and
negative time coincides as a point set with some trajectory of
\(X\) in a neighborhood of \(\partial U\).
\end{Thm}

\medskip

Consider now the semi-space \(\pi^+= \{(x,y,z) \in {\Bbb R}^3 : z
\geq 0\}\) and let \(X: \pi^+ \rightarrow {\Bbb R}^3\) the vector
field \( X(x,y,z)=(-x,-y,z)\). Let \(\phi\) the flow in \(\pi^+\)
induced by \(X\) and let \(D_{m,l}=\{(x,y,z) \in {\Bbb R}^3:
x^2+y^2=1/m, z=1/l\}\).

For every natural number \(n \geq 2\), take the cylinder \(B_n=
\{(x,y,z) \in {\Bbb R}^3 : x^2+y^2=1/n, z \in [1/n, 1/(n-1)]\}\).
Now for every positive even integer \(k\),  we define the flow
boxes \(U_k= \{ \phi(x,t): x \in D_{k,k}, t\geq 0\} \cap B_k\) and
\(V_k = \{ \phi(x,t): x \in D_{k/2,k}, t\geq 0\} \cap B_k\). It is
clear that \(V_k\) is a shrinkage of \(U_k\). On the other hand,
\(U_{k} \cap U_{k'}= \emptyset \) if \(k\ne k'\).

For each \(k \in 2{\Bbb N}\), let \(X_{k}^{1}\) be the vector
field obtained by applying Wilson's theorem to \(X\) and the pair
\((U_k,V_k)\).

Now let \(G:\pi^+ \rightarrow {\Bbb R}^3\) the vector field
defined as \(G(p)=X(p)\) if \(p \notin \bigcup_{k \in 2{\Bbb N}}
U_k\) and \(G(p) = X_{k}^{1}(p)\) if \( p \in U_k\). Finally
consider a flat enough (in \(0\)) smooth non-negative real map
\(\gamma\), depending of \(\|p\|^2\), such that
\(\gamma^{-1}(0)=\{0\}\) to obtain \(X_1=\gamma G\) to be smooth.

Let \(\psi\) the flow in \(\pi^+\) associated to \(X_1\).  The set
of periodic orbits of \(\psi\) is countable. Then we can choose a
positive and decreasing sequence \(t_n \to 0\) such that
\(Fix(\psi(t_n, \cdot)) = Per(\psi(t_n,\cdot))= \{0\}\). Since
each \(D_{k/2,k}\) is a section that captures every orbit in
\(int(\pi^+)\) near \(0\), it is clear that \(0\) is Lyapunov
stable.

Now, we shall apply the same construction of Theorem 1 but we will
paste adequately, in every cone, copies of homeomorphims
conjugated to \(\psi(t_n, \cdot): \pi^+ \rightarrow \pi^+\)
instead of homeomorphisms conjugated to the map \(g_{|\pi^+}\) of
Figure 3 and Example 1.

As in Theorem 1, for every $n=2m-1$ odd we have in each sector
\(E(S_n)\) a finite family of identical  cones \(E(T_{j,m})\), \(j
\in \{1, 2, \dots, q_{m}^{c_m}\}\). For every \(m\) there is a
canonical cone \(E(T_m) \subset int(\pi^+)\) which is isometric to
every \(E(T_{j,m})\). Let \(h_{m}: \pi^+ \rightarrow E(T_m)\) be a
homeomorphisms such that for every for every \(\bar x \in
\partial(E(T_m))\),  \(\|h_{m}^{-1}(\bar x)\|=\|\bar x\|\).

Now define the homeomorphisms \(\psi'_m=h_{m} \circ \psi(t_m,
\cdot) \circ h_{m}^{-1}: E(T_m) \rightarrow E(T_m)\).

Begin with a \({\Bbb R}^3\)-homeomorphism (dynamically equivalent
to the homeomorphism \(f_0\) of Theorem 1 in  \({\Bbb R}^3
\setminus \bigcup int(E(T_{j,m}))\)), \(h_1: {\Bbb R}^3
\rightarrow {\Bbb R}^3\) such that \(Fix(h_1)=\{z=0\}\), \(h_1\)
is decreasing in each ray \(\{\lambda \bar x: \lambda \geq 0, \bar
x \in {\Bbb R}^3 \setminus \{z=0\}\}\) and \(h_1\) behaves in each
ray in \(\partial(E(T_{j,m}))\) as \(\psi'_m\).

Replacing, in each cone \(E(T_{j,m})\), \(h_1\) by copies of
\(\psi'_m\) we obtain a \({\Bbb R}^3\)-homeomorphism \(h_0\).

Let \(h= g_0 \circ h_0\). We obtain in this way a \({\Bbb
R}^3\)-homeomorphism such that \(Fix(h)=Per(h)=\{0\}\) and \(0\)
is Lyapunov stable. It is easy to see that also \(h\) is  limit of
a sequence of homeomorphisms for which every closed ball centered
in \(0\) and large enough radius is still an isolating block with
the same exit sets and the same behavior than in Theorem 1. Then,
the sequence of fixed point indices of the iterates of \(h\) and
\(f\) coincide. \(\square\).

\medskip

\medskip

\noindent {\bf Final Remarks.}

\medskip

i) In Theorem 1, if \(B \subset {\Bbb R}^3\) is any closed ball
centered in the origin,  \(Inv(B,f)\) is the closed 2-disc \(B
\cap \{z=0\}\). For this kind of nice compacta there is a
3-dimensional Carath\'eodory's compactification (see \cite{BJ})
and one could try to apply the ideas of Le Calvez to reduce the
problem of the computation of the indices to the case where the
fixed point is an isolated invariant set. Unfortunately this
method is not longer valid because, in this case, the two
associated fixed prime ends are not isolated invariant sets.

In Theorem 2, if \(B \subset {\Bbb R}^3\) is any closed ball
centered in the origin,  \(Inv(B,f)\) contains the union of the
closed 2-disc \(B \cap \{z=0\}\) and a countable family of
circles.

\medskip
\medskip

ii) Consider the restriction to \(\pi^+\) of the homeomorphisms
\(f\) and \(h\) of Theorems 1 and 2. We can define, by symmetry,
 global \({\Bbb R}^3\)-homeomorphisms, \(F\) and \(H\). Now let
 \(S\) the symmetry with respect to the plane \(\pi=\{z=0\}\).
 Now, \(S \circ F\) and \(S \circ H\) are orientation reversing
 homeomorphism such that \(Fix(S \circ F)=Per(S \circ F)=Fix(S \circ
 H)=Per(S \circ H) =\{0\}\).  Of course \(0\) is again Lyapunov
 stable for \(S \circ H\) and, in this case, \(i((S \circ
 F)^{2k+1},0)=i((S \circ
 H)^{2k+1},0)=1\) for every \(k\in {\Bbb N}\). For even iterates
 we have that
 \(i((S \circ F)^{2k},0)=i((S \circ H)^{2k},0)=-1+ 2i(f^{2k},0)\).

\medskip
\medskip

Francisco R. Ruiz del Portal

Departamento de Geometría y Topología, Facultad de
CC.Mate\-m\'{a}\-ti\-cas, Universidad Complutense de Madrid,
Madrid 28040, Spain.

E-mail: R\(_{-}\)Portal@mat.ucm.es

\medskip
\medskip

José Manuel Salazar.

Departamento de Matemáticas. Universidad de Alcalá.  Alcalá
de Henares. Madrid 28871, Spain.

E-mail: josem.salazar@uah.es

\end{document}